\pdfoutput=1
\documentclass[12pt,a4paper,oneside]{article}
\usepackage{amssymb,latexsym,amsmath,amsthm}
\usepackage{mathtools}
\usepackage[T1]{fontenc}
\usepackage{cmap}
\usepackage[utf8]{inputenc}
\usepackage[english]{babel}

\usepackage{amsfonts}

\usepackage[colorlinks=true,linkcolor=blue,citecolor=blue,urlcolor=blue]{hyperref}

\newcommand{\nid}{\noindent}

\newtheorem{Assertion}{Assertion}
\newtheorem{Theorem}{Theorem}
\newtheorem{Lemma}{Lemma}
\theoremstyle{definition}
\newtheorem{Definition}{Definition}

\begin{document}

\begingroup
\renewcommand{\thefootnote}{}
\footnotetext{\copyright D.~L.~Stupin, V.~G.~Sheretov 2005}
\endgroup
\nid UDC 517.53

\medskip
\begin{center}
\textbf{PROOF OF THE LOCAL KRZYZ CONJECTURE}

\medskip
{D.~L.~Stupin, V.~G.~Sheretov}

{Department of Mathematical Analysis}
\end{center}

\medskip
\begingroup\slshape
\nid For each $t>0$, a special topology is introduced on the class $B_t$ of bounded holomorphic functions $f$
that do not vanish in the unit disk and satisfy $f(0)=e^{-t}$.
This topology is at least as strong as the topology of locally uniform convergence.
A modified method of interior variations is used to prove that, for every index $p$, the functions
$F(e^{i\varphi}z^p,t)$ yield a strict local maximum of the modulus of the $p$-th Taylor coefficient in this topology.
The largest of these local maximum values is $2/e$ and is attained at $t=1$.

This article is a translation of the paper by D.~L.~Stupin and V.~G.~Sheretov,
\emph{Proof of the Local Krzyz Conjecture}, Vestnik TvGU, Series: Applied Mathematics, Tver, 2005,
no.~6(12), pp.~122--125.\par
\endgroup

\section*{Introduction}

Let $\Delta \coloneqq \{z\in\mathbb{C}\,:\,|z|<1\}$ denote the unit disk.
We denote the Taylor coefficients of a function $f$ by $\{f\}_n$, $n\in \{0\}\cup \mathbb N$.

The class $B$ consists of functions $f$ that are holomorphic in $\Delta$, do not vanish there, and satisfy
$|f(z)|\leqslant1$, $z\in\Delta$.

In 1968, the Polish mathematician Jan Krzyz conjectured~\cite{1,2}
that if $f\in B$, then
\[
    |\{f\}_n|\leqslant2/e,
    \qquad
    n\in\mathbb N,
\]
with equality attained by functions of the form
$e^{i\psi}F(e^{i\varphi}z^n,1)$, where
\[
    F(z,t) \coloneqq \exp\left({-t\frac{1+z}{1-z}}\right),
    \qquad
    \varphi,\psi\in\mathbb R,
    \quad
    t\in[0,+\infty).
\]

At present, the Krzyz conjecture has been proved only through the fifth Taylor
coefficient~\cite{3}.

That extremal functions exist for the Krzyz problem at a fixed index $n$ is clear.
Adding the zero function to $B$, namely $f(z)\equiv0$, yields a family
$\overline B$ that is compact in the topology of locally uniform convergence.
The functional assigning to each function $f$ in $\overline B$ the modulus of its
$n$-th Taylor coefficient is continuous on $\overline B$ and therefore
attains its supremum.
If $f(z)=F(-z^n,1)$, then $\{f\}_n=2/e$.

Since the class $B$ is invariant under rotations in the plane
of the variable $w$, we may restrict attention, without loss of generality, to functions with $f(0)>0$.
Since $0<\{f\}_0\leqslant1$, we may set $\{f\}_0=e^{-t}$, with $t\in[0,+\infty)$.
We denote these subclasses by $B_t$.

We denote the set of functions $h$ holomorphic in the disk $\Delta$, normalized by $h(0)=1$,
and having positive real part in $\Delta$
by $C$ and call it the Caratheodory class.

Every function in $B_t$ can be represented in the form
\[
    f(z)=e^{-t\,h(z)},
    \qquad
    h\in C.
\]
For each $t>0$, this formula establishes a bijection between the classes $C$ and $B_t$.
We denote this bijection by $G_t$.

The class $B_0$ consists of the single function $f\equiv1$, so $B_0$ can be regarded as completely
understood.

\section{Preliminaries on Topology}

Let $X$ be a set on which we wish to define a topology.

\begin{Definition}
    Suppose that to each point $x\in X$ we assign a nonempty family $\mathcal N(x)$ of subsets of $X$.
    We call a set $G\subseteq X$ open if for every $x\in G$ there is a set $A\in\mathcal N(x)$
    such that $A\subseteq G$.
\end{Definition}

\begin{Definition}
    A set $U\subseteq X$ is called a neighborhood of a point $x\in X$ if there is an open set $G$
    such that $x\in G\subseteq U$.
\end{Definition}

\begin{Definition}
    A family $\mathcal N(x)$ of neighborhoods of a point $x\in X$ is called a neighborhood base at $x$
    if for every neighborhood $U$ of $x$ there is a set $A\in\mathcal N(x)$ such that
    $A\subseteq U$.
\end{Definition}

\begin{Definition}
    A family $\{\mathcal N(x)\}_{x\in X}$ is called a neighborhood system on $X$
    if $\mathcal N(x)$ is a neighborhood base at $x$ for every $x\in X$.
\end{Definition}

\smallskip
To define a topology using a neighborhood system, we use the following assertion
(see~\cite[p.~25]{Kantorovich:ref1}):

\begin{Assertion}
    \label{st:Stupin:neighborhoodSystemTopology}
    Let $X$ be a set, and suppose that to every point $x\in X$ there corresponds a nonempty family $\mathcal N(x)$
    of subsets of $X$ satisfying the following conditions\textnormal:
    \newline\textnormal1.
    $x\in A$ for every $A\in\mathcal N(x)$\textnormal;
    \newline\textnormal2.
    if $A_1,A_2\in\mathcal N(x)$, then there exists a set $A_3\in\mathcal N(x)$
    such that $A_3\subseteq A_1\cap A_2$\textnormal;
    \newline\textnormal3.
    if $A_x\in\mathcal N(x)$, then there exists a set $A_x'\in\mathcal N(x)$ such that $A_x'\subseteq A_x$ and
    for every point $y\in A_x'$ there exists a set $A_y\in\mathcal N(y)$ such that $A_y\subseteq A_x$.
    \newline
    Then the collection of open subsets of $X$ is a topology on $X$
    for which $\{\mathcal N(x)\}_{x\in X}$ is a neighborhood system.
\end{Assertion}

\section{The Special Topology}

Fix $t>0$ and define a topology on $B_t$.
We first define a topology on $C$ and then transfer it to $B_t$ via the bijection $G_t$.

For $H,h\in C$ and $\varepsilon_0\in(0,1]$, set
\[
    l_{\varepsilon_0}(H,h) \coloneqq \{(1-\varepsilon)H+\varepsilon h
    \,\colon\,\varepsilon\in[0,\varepsilon_0)\},
    \qquad
    U_{\varepsilon_0}(H) \coloneqq \bigcup\limits_{h\in C}l_{\varepsilon_0}(H,h).
\]
By convexity of the class $C$, we have $U_{\varepsilon_0}(H)\subset C$.
Set $\mathcal N_C(H)\coloneqq \{U_{\varepsilon_0}(H):0<\varepsilon_0\leqslant1\}$.

\begin{Assertion}
The family $\{\mathcal N_C(H)\}_{H\in C}$ generates a topology on $C$ and is a neighborhood system
for that topology.
\end{Assertion}

\begin{proof}
We verify that $\{\mathcal N_C(H)\}_{H\in C}$ satisfies all the conditions of
Assertion~\ref{st:Stupin:neighborhoodSystemTopology}.

\noindent
\textbf{1.}
The first condition holds because $H\in U_{\varepsilon_0}(H)$.

\noindent
\textbf{2.}
The second condition holds because
$U_{\min\{\varepsilon_1,\varepsilon_2\}}(H)\subseteq U_{\varepsilon_1}(H)\cap U_{\varepsilon_2}(H)$.

\noindent
\textbf{3.}
To verify the third condition, it suffices to show that for every point $q\in U_{\varepsilon_0}(H)$ there exists
a number $\mu\in(0,1]$ such that $U_\mu(q)\subseteq U_{\varepsilon_0}(H)$.
Since $q\in U_{\varepsilon_0}(H)$, we have
$q
=
(1-\alpha)H+\alpha g$,
$g\in C$,
$0\leqslant\alpha<\varepsilon_0$.
We show that we can take $\mu\coloneqq\tfrac{\varepsilon_0-\alpha}{1-\alpha}$.
Since $0\leqslant\alpha<\varepsilon_0\leqslant1$, we have $0<\mu\leqslant\varepsilon_0\leqslant1$.

By definition, $U_\mu(q)$ consists of functions of the form
\[
    u
    =
    (1-\beta)q+\beta p
    =
    (1-\beta)(1-\alpha)H+(1-\beta)\alpha g+\beta p,
    \quad
    p\in C,
    \quad
    0\leqslant\beta<\mu\leqslant1.
\]

We show that $u=(1-\gamma)H+\gamma w$ with $w\in C$.
Let $\gamma\coloneqq(1-\beta)\alpha+\beta$ be the sum of the coefficients of $g$ and $p$.
Clearly, $0\leqslant\gamma<1$.
If $\gamma=0$, then $\alpha=\beta=0$ and $u=H$.
We may therefore set $w\equiv1$.
If $\gamma>0$, set
$w
\coloneqq
\frac{(1-\beta)\alpha}{\gamma}g
+
\frac{\beta}{\gamma}p$.
The coefficients of $g$ and $p$ are nonnegative and sum to $1$, so $w\in C$ by convexity of the class $C$.

By definition, $u\in U_{\varepsilon_0}(H)$ if $0\leqslant\gamma<\varepsilon_0$.
If $\gamma=0$, this inequality plainly holds. If $\gamma>0$, it holds for
$\beta<\mu$.
If $u\in U_\mu(q)$, then $\beta<\mu$, so $\gamma<\varepsilon_0$ and $u\in U_{\varepsilon_0}(H)$.
Consequently, $U_\mu(q)\subseteq U_{\varepsilon_0}(H)$.

In Assertion~\ref{st:Stupin:neighborhoodSystemTopology}, we may take $A_H=A_H'=U_{\varepsilon_0}(H)$,
so its third condition holds.

Thus all the conditions of Assertion~\ref{st:Stupin:neighborhoodSystemTopology} hold,
which proves the claim.
\end{proof}

\smallskip
We denote the resulting topology by $\tau_C$.
We have $G_t(\emptyset)=\emptyset$ and $G_t(C)=B_t$.
The image of a union of subsets of $C$ is the union of their images.
Since $G_t$ is bijective, the image of the intersection of two subsets of $C$ is the intersection of their images.
Therefore the collection
\[
    \tau_{B_t} \coloneqq \{G_t(U_{\alpha})\,:\,U_{\alpha}\in\tau_C\}
\]
is a topology on the class $B_t$.

\begin{Assertion}
The topology $\tau_{B_t}$ is at least as strong as the topology of locally uniform convergence in $\Delta$
induced on $B_t$.
\end{Assertion}

\begin{proof}
Fix a function $f=G_t(H)\in B_t$, a nonempty compact set $K\subset\Delta$, and a number $\xi>0$.
Since the set
$\big\{\widetilde f\in B_t\,\colon\,\max\limits_{z\in K}|\widetilde f(z)-f(z)|<\xi\big\}$
is a basic neighborhood of $f$ in the topology of locally uniform convergence
induced on $B_t$, it suffices to find $\varepsilon_0\in(0,1]$ such that for every function
$\widetilde f\in G_t(U_{\varepsilon_0}(H))$ the inequality
$\max\limits_{z\in K}|\widetilde f(z)-f(z)|<\xi$.

Let $r\coloneqq\max\limits_{z\in K}|z|$ and $M_K\coloneqq\frac{1+r}{1-r}$.
For any function $h\in C$, Schwarz's lemma gives $\big|\frac{h(z)-1}{h(z)+1}\big|\leqslant|z|$.
It follows that $\max\limits_{z\in K}|h(z)|\leqslant M_K$.

Let $\widetilde f=G_t(q)$, where $q=(1-\varepsilon)H+\varepsilon h$,
$h\in C$, $0\leqslant\varepsilon<\varepsilon_0$.
Then $\max\limits_{z\in K}|q(z)-H(z)|<2\varepsilon_0 M_K$.

For each point $z\in K$, the numbers $a=tq(z)$ and $b=tH(z)$ have positive real parts.
The identity
\[
    e^{-a}-e^{-b}=(b-a)\int\limits_0^1 e^{-((1-s)b+sa)}\,ds
\]
implies that $|\widetilde f(z)-f(z)|\leqslant t|q(z)-H(z)|$.
Thus $\max\limits_{z\in K}|\widetilde f(z)-f(z)|<2t\varepsilon_0 M_K$.
Choose $\varepsilon_0\in(0,1]$ so that $2t\varepsilon_0 M_K<\xi$.

Consequently, every neighborhood of $f$ in the topology of locally uniform convergence contains a neighborhood
in the topology $\tau_{B_t}$, which proves the assertion.
\end{proof}

\section{Auxiliary Results}

\begin{Lemma}
Let $q\in\mathbb N$, $\sum\limits_{k=1}^q\alpha_k=1$, and $\alpha_k\geqslant0$ for $k=1,\ldots,q$.
Then
\begin{equation}
    \label{eq:Stupin:alphaSquaresLowerBound}
    \sum\limits_{k=1}^q\alpha_k^2\geqslant\frac1q,
\end{equation}
with equality if and only if
\begin{equation}
    \label{eq:Stupin:alphaEquals}
    \alpha_1=\alpha_2=\ldots=\alpha_q=\frac1q.
\end{equation}
\end{Lemma}

\begin{proof}
Let $\boldsymbol x=(\alpha_1,\ldots,\alpha_q)$ and $\boldsymbol y=(1,\ldots,1)$ be
vectors in $\mathbb R^q$.
Applying the Cauchy--Schwarz inequality to the vectors $\boldsymbol x$ and $\boldsymbol y$ gives
\[
    1
    =
    \left(\sum\limits_{k=1}^q\alpha_k\right)^2
    =
    \left(\sum\limits_{k=1}^q\alpha_k\cdot1\right)^2
    \leqslant
    \sum\limits_{k=1}^q\alpha_k^2
    \cdot
    \sum\limits_{k=1}^q1^2
    =
    q\sum\limits_{k=1}^q\alpha_k^2.
\]
This yields inequality~\eqref{eq:Stupin:alphaSquaresLowerBound}.
Equality in the Cauchy--Schwarz inequality holds if and only if
the vectors $\boldsymbol x$ and $\boldsymbol y$ are collinear, so
equality in~\eqref{eq:Stupin:alphaSquaresLowerBound} holds in precisely those cases,
that is, if and only if $\alpha_1=\ldots=\alpha_q$.
These equalities and the conditions $\sum\limits_{k=1}^q\alpha_k=1$ and $\alpha_k\geqslant0$ for $k=1,\ldots,q,$
give~\eqref{eq:Stupin:alphaEquals}.
\end{proof}

\begin{Lemma}
\label{lem:Stupin:caratheodoryCoeffDifference}
Let $h\in C$, $p\in\mathbb N$, and $1\leqslant j<p$.
Then
\[
    |\{h\}_{p-j}-\overline{\{h\}_j}|
    \leqslant 2\sqrt{2-\operatorname{Re}\{h\}_p}.
\]
\end{Lemma}

\begin{proof}
Let $0<r<1$.
By Cauchy's formula,
\begin{multline*}
    r^{p-j}\{h\}_{p-j}-r^j\overline{\{h\}_j}
    =
    \frac1{2\pi}\int\limits_0^{2\pi}e^{-i(p-j)\theta}h(re^{i\theta})\,d\theta
    -
    \frac1{2\pi}\int\limits_0^{2\pi}e^{ij\theta}\overline{h(re^{i\theta})}\,d\theta
    =\\=
    \frac1{2\pi}\int\limits_0^{2\pi}
    \bigl(e^{-i(p-j)\theta}-e^{ij\theta}\bigr)
    \bigl(h(re^{i\theta})+\overline{h(re^{i\theta})}\bigr)\,d\theta
    -\\-
    \frac1{2\pi}\overline{\int\limits_0^{2\pi}e^{i(p-j)\theta}h(re^{i\theta})\,d\theta}
    +
    \frac1{2\pi}\int\limits_0^{2\pi}e^{ij\theta}h(re^{i\theta})\,d\theta
    =\\=
    \frac1\pi\int\limits_0^{2\pi}
    e^{ij\theta}(e^{-ip\theta}-1)\operatorname{Re}h(re^{i\theta})\,d\theta,
\end{multline*}
since
$\int\limits_0^{2\pi}e^{i(p-j)\theta}h(re^{i\theta})\,d\theta=
\int\limits_0^{2\pi}e^{ij\theta}h(re^{i\theta})\,d\theta=0$
by the mean value theorem for the functions $z^{p-j}h(z)$ and $z^j h(z)$, which vanish at $z=0$.

Let $u(\theta)=|e^{-ip\theta}-1|\sqrt{\frac1{2\pi}\operatorname{Re}h(re^{i\theta})}$ and
$v(\theta)=\sqrt{\frac1{2\pi}\operatorname{Re}h(re^{i\theta})}$ be continuous functions.
Applying the Cauchy--Schwarz inequality to $u$ and $v$ and noting that
$\frac1{2\pi}\int\limits_0^{2\pi}\operatorname{Re}h(re^{i\theta})\,d\theta=\operatorname{Re}h(0)=1$
by the mean value theorem, we obtain
\begin{multline*}
    \big|r^{p-j}\{h\}_{p-j}-r^j\overline{\{h\}_j}\big|
    =\\=
    \left|\frac1\pi\int\limits_0^{2\pi}
    e^{ij\theta}(e^{-ip\theta}-1)\operatorname{Re}h(re^{i\theta})\,d\theta\right|
    \leqslant
    \frac1{\pi}\int\limits_0^{2\pi}
    |e^{-ip\theta}-1|\operatorname{Re}h(re^{i\theta})\,d\theta
    \leqslant\\\leqslant
    2
    \cdot
    \sqrt{\frac1{2\pi}\int\limits_0^{2\pi}|e^{-ip\theta}-1|^2
    \operatorname{Re}h(re^{i\theta})\,d\theta}
    \cdot
    \sqrt{\frac1{2\pi}\int\limits_0^{2\pi}\operatorname{Re}h(re^{i\theta})\,d\theta}
    =\\=
    2
    \cdot
    \sqrt{\frac1{2\pi}\int\limits_0^{2\pi}|e^{-ip\theta}-1|^2
    \operatorname{Re}h(re^{i\theta})\,d\theta}.
\end{multline*}

Since $\frac1{2\pi}\int\limits_0^{2\pi}e^{-ip\theta}h(re^{i\theta})\,d\theta=r^p\{h\}_p$ by Cauchy's formula
and $\int\limits_0^{2\pi}e^{ip\theta}h(re^{i\theta})\,d\theta=0$ by the mean value theorem for
the function $z^p h(z)$, we have

\begin{multline*}
    \frac1{2\pi}\int\limits_0^{2\pi}|e^{-ip\theta}-1|^2
    \operatorname{Re}h(re^{i\theta})\,d\theta
    =\\=
    \frac1{2\pi}\int\limits_0^{2\pi}(2-2\operatorname{Re}e^{-ip\theta})
    \operatorname{Re}h(re^{i\theta})\,d\theta
    =
    2-2\operatorname{Re}\left(\frac1{2\pi}\int\limits_0^{2\pi}
    e^{-ip\theta}\operatorname{Re}h(re^{i\theta})\,d\theta\right)
    =\\=
    2-\operatorname{Re}\left(
    \frac1{2\pi}\int\limits_0^{2\pi}e^{-ip\theta}h(re^{i\theta})\,d\theta
    +
    \frac1{2\pi}\overline{\int\limits_0^{2\pi}e^{ip\theta}h(re^{i\theta})\,d\theta}
    \right)
    =
    2-r^p\operatorname{Re}\{h\}_p.
\end{multline*}

Letting $r$ tend to $1$ gives the desired inequality.
\end{proof}

For $p\in\mathbb N$, set $H_p(z)\coloneqq\frac{1+z^p}{1-z^p}$.

If $f$ and $g$ are functions holomorphic in $\Delta$, then
\begin{equation}
    \label{eq:Stupin:prodCoeff}
    \{f \cdot g\}_p
    =
    \{f\}_p\{g\}_0+\{f\}_{p-1}\{g\}_1+\ldots+\{f\}_0\{g\}_p.
\end{equation}

\begin{Lemma}
\label{lem:Stupin:generalVariation}
Let $t>0$, $p\in\mathbb N$, $h\in C$, and $0<\varepsilon<1$.
Then
\begin{equation}
    \label{eq:Stupin:generalVariation}
    \{e^{-t\varepsilon(h-H_p)}\}_p
    =
    t(2-\{h\}_p)\varepsilon
    +\frac{t^2}{2}\left(\sum\limits_{j=1}^{p-1}\{h\}_j\{h\}_{p-j}\right)\varepsilon^2
    +R_{p,\varepsilon}(h),
\end{equation}
where $|R_{p,\varepsilon}(h)|\leqslant C_{t,p}\varepsilon^3\sum\limits_{j=1}^{p-1}|\{h\}_j|^2$,
and the constant $C_{t,p}\geqslant0$ depends only on $t$ and $p$.
\end{Lemma}

\begin{proof}
We have
\[
    e^{-t\varepsilon(h(z)-H_p(z))}
    =
    \sum\limits_{m=0}^{\infty}\frac{(-t)^m}{m!}(h(z)-H_p(z))^m\varepsilon^m.
\]
It follows that
\begin{equation*}
    \{e^{-t\varepsilon(h-H_p)}\}_p
    =
    -t\{h-H_p\}_p\varepsilon
    +\frac{t^2}{2}\{(h-H_p)^2\}_p\varepsilon^2
    +
    \sum\limits_{m=3}^{p}\frac{(-t)^m}{m!}\{(h-H_p)^m\}_p\varepsilon^m.
\end{equation*}

We now calculate the terms in~\eqref{eq:Stupin:generalVariation} in order.

\textbf{1.}
We have $\{h-H_p\}_p=\{h\}_p-2$.

\textbf{2.}
Since $\{h-H_p\}_0=0$ and $\{H_p\}_j=0$ for $j=1,\ldots,p-1$,
formula~\eqref{eq:Stupin:prodCoeff} gives
\[
    \{(h-H_p)^2\}_p
    =
    \sum\limits_{j=0}^p\{h-H_p\}_j\{h-H_p\}_{p-j}
    =
    \sum\limits_{j=1}^{p-1}\{h\}_j\{h\}_{p-j}.
\]

\textbf{3.}
Since $(h-H_p)^m$ is a product of $m$ factors $h-H_p$,
applying~\eqref{eq:Stupin:prodCoeff} $m-1$ times for $3\leqslant m\leqslant p$ yields
\[
    \{(h-H_p)^m\}_p
    =
    \sum\limits_{\substack{j_1+\cdots+j_m=p\\j_1,\ldots,j_m\geqslant1}}
    \{h-H_p\}_{j_1}\cdots\{h-H_p\}_{j_m}
    =
    \sum\limits_{\substack{j_1+\cdots+j_m=p\\j_1,\ldots,j_m\geqslant1}}
    \{h\}_{j_1}\cdots\{h\}_{j_m}.
\]
The indices $j_i$ start at one since $\{h-H_p\}_0=0$.
Moreover, each $j_i<p$, so $\{h-H_p\}_{j_i}=\{h\}_{j_i}$.

The Caratheodory estimate (see~\cite[p.~57]{5}) gives $|\{h\}_j|\leqslant2$.
For $a,b\in\mathbb R$, the inequality $2ab\leqslant a^2+b^2$ follows from $(a-b)^2\geqslant0$.
Applying it with $a=|\{h\}_{j_1}|$ and $b=|\{h\}_{j_2}|$ gives
\begin{equation*}
    |\{h\}_{j_1}\cdots\{h\}_{j_m}|
    \leqslant 2^{m-2}|\{h\}_{j_1}\{h\}_{j_2}|
    \leqslant
    2^{m-2}\frac{|\{h\}_{j_1}|^2+|\{h\}_{j_2}|^2}{2}
    \leqslant 2^{m-2}\sum\limits_{j=1}^{p-1}|\{h\}_j|^2.
\end{equation*}
Consequently,
\begin{multline*}
    |R_{p,\varepsilon}(h)|
    =
    \left|\sum\limits_{m=3}^{p}\frac{(-t)^m}{m!}\{(h-H_p)^m\}_p\varepsilon^m\right|
    \leqslant
    \sum\limits_{m=3}^{p}\frac{t^m}{m!}
    \sum\limits_{\substack{j_1+\cdots+j_m=p\\j_1,\ldots,j_m\geqslant1}}
    |\{h\}_{j_1}\cdots\{h\}_{j_m}|\,\varepsilon^m
    \leqslant\\\leqslant
    C_{t,p}\varepsilon^3\sum\limits_{j=1}^{p-1}|\{h\}_j|^2,
    \quad\text{where}\quad
    C_{t,p}
    \coloneqq
    \sum\limits_{m=3}^{p}\frac{2^{m-2}t^m}{m!}
        \sum\limits_{\substack{j_1+\cdots+j_m=p\\j_1,\ldots,j_m\geqslant1}}1.
\end{multline*}
Note that $R_{1,\varepsilon}(h)=0$ and $R_{2,\varepsilon}(h)=0$.
\end{proof}

\section{Main Result}

In~\cite{4} it was proved that $F(-z^p,1)$ gives a strict local maximum of the functional
$\operatorname{Re}\{f\}_p$ on $\bigcup\limits_{t\in[0,+\infty)}B_t$
in the topology of locally uniform convergence.

We prove an analogous statement for the modulus of the $p$-th coefficient
in the topology $\tau_{B_t}$ on each class $B_t$ with $t>0$.

\begin{Theorem}
Let $t>0$ and $p\in\mathbb N$.
Then the functions $F(e^{i\varphi}z^p,t)$, $\varphi\in[0,2\pi)$,
give a strict local maximum of the functional $I_{t,p}[f]\coloneqq|\{f\}_p|$ on the class $B_t$
in the topology $\tau_{B_t}$.
Moreover, $I_{t,p}[F(e^{i\varphi}z^p,t)]=2te^{-t}$, and this value attains its maximum of $2/e$ at $t=1$.
\end{Theorem}

\begin{proof}
\textbf{1.}
Let $t>0$, $p\in\mathbb N$, and
\[
    H_p(z) \coloneqq \frac{1+z^p}{1-z^p},
    \qquad
    f_{t,p}(z) \coloneqq e^{-tH_p(z)}.
\]
We show that the function $f_{t,p}$ gives a strict local maximum of the functional $I_{t,p}[f]$
on $B_t$ in the topology $\tau_{B_t}$.
To this end, consider the variation
\[
    f_{t,p,\varepsilon,h}(z)
    \coloneqq
    e^{-t((1-\varepsilon)H_p(z)+\varepsilon h(z))},
    \qquad
    h\in C,
    \quad
    \varepsilon\in(0,1)
\]
of the function $f_{t,p}$.
It follows that
\[
    \{f_{t,p,\varepsilon,h}\}_p=-\frac{t}{e^t}
    (2-(2-\{h\}_p)\varepsilon+o(\varepsilon))
\]
as $\varepsilon\to0$.

\textbf{2.}
Let $h\in C$ and $\{h\}_p\neq2$.
Writing $w\coloneqq-(2-\{h\}_p)\varepsilon+o(\varepsilon)$, we have
\begin{multline*}
    |2+w|
    =
    \sqrt{(2+w)(2+\overline w)}
    =
    \sqrt{4+4\operatorname{Re}w+|w|^2}
    =\\=
    \sqrt{4+4\operatorname{Re}w+o(\varepsilon)}
    =
    2+\operatorname{Re}w+o(\varepsilon).
\end{multline*}
That is, $|2-(2-\{h\}_p)\varepsilon+o(\varepsilon)|=2-(2-\operatorname{Re}\{h\}_p)\varepsilon+o(\varepsilon)$.

Therefore, for sufficiently small $\varepsilon>0$,
\[
    |\{f_{t,p,\varepsilon,h}\}_p|
    =
    \frac{t}{e^t}(2-(2-\operatorname{Re}\{h\}_p)\varepsilon+o(\varepsilon))
    <
    |\{f_{t,p}\}_p|.
\]
Indeed, since $|\{h\}_p|\leqslant2$ and $\{h\}_p\neq2$, we have $2-\operatorname{Re}\{h\}_p>0$.
Thus, for sufficiently small $\varepsilon>0$, this variation decreases $|\{f_{t,p}\}_p|$.

\textbf{3.}
We call the variation $f_{t,p,\varepsilon,h}$ trivial if $f_{t,p,\varepsilon,h}=f_{t,p}$.

Now let $h\in C$ and $\{h\}_p=2$.
According to~\cite[p.~57]{5}, this is possible only when
\[
    h(z)
    =
    h_q(z)
    \coloneqq
    \sum\limits_{k=1}^q\alpha_k
    \frac{1+e^{\frac{2\pi ki}{q}}z}{1-e^{\frac{2\pi ki}{q}}z},
\]
where $\sum\limits_{k=1}^q\alpha_k=1$, \ $\alpha_k\geqslant0$, \
$k=1,\ldots,q$, \ and $p$ is a multiple of $q$.
We vary the function $f_{t,p}$ using $h_q$.
If $q=p$ and
\[
    \alpha_1=\alpha_2=\ldots=\alpha_q=\frac1q,
\]
the variation is trivial, so we exclude this case.
We have
\begin{multline*}
    f_{t,p,\varepsilon,h}(z)
    \coloneqq
    \exp({-t((1-\varepsilon)H_p(z)+\varepsilon h_q(z))})
    =\\=
    e^{-t}\exp\left({-2t\left(\sum\limits_{j=1}^{\infty}z^{jp}+
    \varepsilon\sum\limits_{\substack{j=1\\ j\not\div p}}^{\infty}
    \left(\sum\limits_{k=1}^q\alpha_k e^{j\frac{2\pi ki}{q}}\right)z^j\right)}\right)
    =\\=
    e^{-t}\exp\left({-2t\left(z^p+\ldots+\varepsilon\left(
    \sum\limits_{j=1}^{p-1}\left(\sum\limits_{k=1}^q
    \alpha_k e^{j\frac{2\pi ki}{q}}\right)z^j+\ldots\right)\right)}\right).
\end{multline*}
For this variation of $f_{t,p}$, the $p$-th Taylor coefficient
is given by
\begin{multline*}
    \{f_{t,p,\varepsilon,h}\}_p
    =
    -\frac{t}{e^t}
    \left[2-2t\left(\sum\limits_{\substack{j,l=1\\ j+l=p}}^{p-1}
    \left(\sum\limits_{k=1}^q\alpha_k e^{j\frac{2\pi ki}{q}}\right)
    \left(\sum\limits_{m=1}^q\alpha_m e^{l\frac{2\pi mi}{q}}\right)
    \right)\varepsilon^2+o(\varepsilon^2)\right]
    =\\=
    -\frac{t}{e^t}\left[2-2t\left(\sum\limits_{j=1}^{p-1}
    \left(\sum\limits_{k=1}^q\alpha_k e^{j\frac{2\pi ki}{q}}\right)
    \left(\sum\limits_{m=1}^q\alpha_m e^{(p-j)\frac{2\pi mi}{q}}\right)\right)
    \varepsilon^2+o(\varepsilon^2)\right]
    =\\=
    -\frac{t}{e^t}\left[2-2t\left(\sum\limits_{j=1}^{p-1}
    \sum\limits_{k,m=1}^q\alpha_k\alpha_m e^{j\frac{2\pi(k-m)i}q}\right)
    \varepsilon^2+o(\varepsilon^2)\right]
\end{multline*}
as $\varepsilon\to0$.

We rewrite the double sum in parentheses.
Since $p$ is a multiple of $q$, we have
\begin{multline*}
    \sum\limits_{j=1}^{p-1}
    \sum\limits_{k,m=1}^q\alpha_k\alpha_me^{j\frac{2\pi(k-m)i}q}
    =\\=
    \sum\limits_{j=1}^{p-1}\left(\sum\limits_{k=1}^q\alpha_k^2+
    \sum\limits_{\substack{k,m=1\\ k<m}}^q\alpha_k\alpha_m e^{j\frac{2\pi(k-m)i}q}+
    \sum\limits_{\substack{k,m=1\\ m<k}}^q\alpha_k\alpha_m e^{j\frac{2\pi(k-m)i}q}
    \right)
    =\\=
    \sum\limits_{j=1}^{p-1}\left(\sum\limits_{k=1}^q\alpha_k^2+
    \sum\limits_{\substack{k,m=1\\ k<m}}^q\alpha_k\alpha_m
    \left(e^{j\frac{2\pi(k-m)i}q}+e^{j\frac{2\pi(m-k)i}q}\right)\right)
    =\\=
    \sum\limits_{j=1}^{p-1}\left(\sum\limits_{k=1}^q\alpha_k^2+
    2\sum\limits_{\substack{k,m=1\\ k<m}}^q\alpha_k\alpha_m
    \cos\left(j\frac{2\pi(m-k)}q\right)\right)
    =\displaybreak[1]\\=
    (p-1)\sum\limits_{k=1}^q\alpha_k^2+
    2\sum\limits_{\substack{k,m=1\\ k<m}}^q\alpha_k\alpha_m
    \sum\limits_{j=1}^{p-1}\cos\left(j\frac{2\pi(m-k)}q\right)
    =\\=
    (p-1)\sum\limits_{k=1}^q\alpha_k^2+
    2\sum\limits_{\substack{k,m=1\\ k<m}}^q\alpha_k\alpha_m
    \left(\frac{\sin\frac{p(m-k)\pi}{q}\cos\frac{(p-1)(m-k)\pi}{q}}
    {\sin\frac{(m-k)\pi}{q}}-1\right)
    =\\=
    (p-1)\sum\limits_{k=1}^q\alpha_k^2-
    2\sum\limits_{\substack{k,m=1\\ k<m}}^q\alpha_k\alpha_m
    =
    (p-1)\sum\limits_{k=1}^q\alpha_k^2-
    \left(1-\sum\limits_{k=1}^q\alpha_k^2\right)
    =
    p\sum\limits_{k=1}^q\alpha_k^2-1.
\end{multline*}
Thus
\begin{equation}
    \label{eq:Stupin:secondVariation}
    \{f_{t,p,\varepsilon,h}\}_p
    =
    -\frac{t}{e^t}\left[2-2t\left(
    p\sum\limits_{k=1}^q\alpha_k^2-1\right)
    \varepsilon^2+o(\varepsilon^2)\right]
\end{equation}
as $\varepsilon\to0$.

\textbf{4.}
If $q<p$, then by~\eqref{eq:Stupin:alphaSquaresLowerBound},
\[
    p\sum\limits_{k=1}^q\alpha_k^2-1\geqslant\frac{p}{q}-1>0.
\]

If $q=p$ and the variation is nontrivial, condition~\eqref{eq:Stupin:alphaEquals} fails.
Therefore the equality condition in~\eqref{eq:Stupin:alphaSquaresLowerBound} implies that
\[
    \sum\limits_{k=1}^q\alpha_k^2>\frac1q,
\]
whence
\[
    p\sum\limits_{k=1}^q\alpha_k^2-1>\frac{p}{q}-1=0.
\]

Thus in either case
\begin{equation}
    \label{eq:Stupin:secondVariationPositive}
    c\coloneqq 2t\left(p\sum\limits_{k=1}^q\alpha_k^2-1\right)>0.
\end{equation}

\textbf{5.}
By~\eqref{eq:Stupin:secondVariation},
\[
    \{f_{t,p,\varepsilon,h}\}_p=-\frac{t}{e^t}(2-c\varepsilon^2+o(\varepsilon^2)).
\]
For sufficiently small $\varepsilon>0$, we have $2-c\varepsilon^2>0$ and
$|o(\varepsilon^2)|\leqslant\frac{c}2\varepsilon^2$.
Therefore,
\begin{equation*}
    |\{f_{t,p,\varepsilon,h}\}_p|
    =\frac{t}{e^t}|2-c\varepsilon^2+o(\varepsilon^2)|
    \leqslant
    \frac{t}{e^t}(2-c\varepsilon^2+|o(\varepsilon^2)|)
    \leqslant\frac{t}{e^t}\left(2-\frac{c}2\varepsilon^2\right)
    <\frac{2t}{e^t}
    =|\{f_{t,p}\}_p|.
\end{equation*}
Thus, for sufficiently small $\varepsilon>0$, this variation decreases the modulus
of the $p$-th coefficient of $f_{t,p}$.

\textbf{6.}
Suppose that $f_{t,p}$ does not give a strict local maximum of the functional $I_{t,p}$
in the topology $\tau_{B_t}$.
Then for every $n\in\mathbb N$ there exist a function $h\in C$, $h\neq H_p$, and
a function $f_{t,p,\varepsilon',h}\in G_t(U_{1/n}(H_p))$, $f_{t,p,\varepsilon',h}\neq f_{t,p}$, such that
$|\{f_{t,p,\varepsilon',h}\}_p|\geqslant|\{f_{t,p}\}_p|$, $\varepsilon'\in(0,1/n)$.

By steps~2 and~5, for the same $h$ and for all sufficiently small $\varepsilon\in(0,\varepsilon')$,
\[
    |\{f_{t,p,\varepsilon,h}\}_p|<2te^{-t}.
\]
For fixed $h$, the function $\varepsilon\mapsto |\{f_{t,p,\varepsilon,h}\}_p|$ is continuous
on the interval $[0,\varepsilon']$, where $f_{t,p,0,h}\coloneqq f_{t,p}$.
Indeed, $\{f_{t,p,\varepsilon,h}\}_p$ is a polynomial in $\varepsilon$
(its degree is at most $p$ because $h(0)=H_p(0)$).
At $\varepsilon=\varepsilon'$ we have $|\{f_{t,p,\varepsilon',h}\}_p|\geqslant2te^{-t}$.
By the intermediate value theorem, there is therefore an $\varepsilon\in(0,\varepsilon']\subset(0,1/n)$ for which
\begin{equation}
    \label{eq:Stupin:returnValue}
    |\{f_{t,p,\varepsilon,h}\}_p|=2te^{-t}.
\end{equation}
We show that there exists a number $\varepsilon_0>0$ such that for every $h\in C$ with $h\neq H_p$ and
every $0<\varepsilon<\varepsilon_0$, equality~\eqref{eq:Stupin:returnValue} is impossible.

Since $f_{t,p,\varepsilon,h}=f_{t,p}e^{-t\varepsilon(h-H_p)}$,
$f_{t,p}(z)=e^{-t}(1-2tz^p)+O(z^{2p})$ as $z\to0$, and $h(0)=H_p(0)$, the
product formula~\eqref{eq:Stupin:prodCoeff} and Lemma~\ref{lem:Stupin:generalVariation} give
\begin{multline}
    \label{eq:Stupin:returnExpansion}
    -\frac{e^t}{t}\{f_{t,p,\varepsilon,h}\}_p
    =
    2-\frac1t\{e^{-t\varepsilon(h-H_p)}\}_p
    =\\=
    2-(2-\{h\}_p)\varepsilon
    -\frac t2\left(\sum\limits_{j=1}^{p-1}\{h\}_j\{h\}_{p-j}\right)\varepsilon^2
    +O\left(\varepsilon^3\sum\limits_{j=1}^{p-1}|\{h\}_j|^2\right).
\end{multline}

We have
\[
    \left|\sum\limits_{j=1}^{p-1}\{h\}_j\{h\}_{p-j}\right|
    \leqslant\sum\limits_{j=1}^{p-1}|\{h\}_j\{h\}_{p-j}|
    \leqslant\frac12\sum\limits_{j=1}^{p-1}\bigl(|\{h\}_j|^2+|\{h\}_{p-j}|^2\bigr)
    =\sum\limits_{j=1}^{p-1}|\{h\}_j|^2.
\]
It follows from~\eqref{eq:Stupin:returnExpansion} that
\begin{equation}
    \label{eq:Stupin:returnDifference}
    -\frac{e^t}{t}\{f_{t,p,\varepsilon,h}\}_p-\bigl(2-\varepsilon(2-\{h\}_p)\bigr)
    =O\left(\varepsilon^2\sum\limits_{j=1}^{p-1}|\{h\}_j|^2\right).
\end{equation}

Since $|\{h\}_p|\leqslant2$, we have
\[
    |2-\varepsilon(2-\{h\}_p)|
    =|2(1-\varepsilon)+\varepsilon\{h\}_p|
    \leqslant2(1-\varepsilon)+\varepsilon|\{h\}_p|
    \leqslant2.
\]
From~\eqref{eq:Stupin:returnValue} and~\eqref{eq:Stupin:returnDifference} we obtain
\begin{equation*}
    0\leqslant 2-|2-\varepsilon(2-\{h\}_p)|
    \leqslant
    \left|-\frac{e^t}{t}\{f_{t,p,\varepsilon,h}\}_p-\bigl(2-\varepsilon(2-\{h\}_p)\bigr)\right|
    =
    O\left(\varepsilon^2\sum\limits_{j=1}^{p-1}|\{h\}_j|^2\right).
\end{equation*}

Moreover, from $|2-\{h\}_p|^2
=4-4\operatorname{Re}\{h\}_p+|\{h\}_p|^2
\leqslant 8-4\operatorname{Re}\{h\}_p
=4(2-\operatorname{Re}\{h\}_p)$ we obtain
\begin{multline*}
    4\varepsilon(1-\varepsilon)(2-\operatorname{Re}\{h\}_p)
    \leqslant4\varepsilon(2-\operatorname{Re}\{h\}_p)-\varepsilon^2|2-\{h\}_p|^2
    =
    4-|2-\varepsilon(2-\{h\}_p)|^2
    =\\=
    \bigl(2-|2-\varepsilon(2-\{h\}_p)|\bigr)\bigl(2+|2-\varepsilon(2-\{h\}_p)|\bigr)
    \leqslant4\bigl(2-|2-\varepsilon(2-\{h\}_p)|\bigr).
\end{multline*}
Consequently, if~\eqref{eq:Stupin:returnValue} holds and $0<\varepsilon<1/2$, then
\begin{equation}
    \label{eq:Stupin:returnDeficit}
    2-\operatorname{Re}\{h\}_p
    =O\left(\varepsilon\sum\limits_{j=1}^{p-1}|\{h\}_j|^2\right).
\end{equation}

Since $|\{h\}_{p-j}-\overline{\{h\}_j}|^2
=
|\{h\}_{p-j}|^2+|\{h\}_j|^2-2\operatorname{Re}(\{h\}_j\{h\}_{p-j})$, summing
over $j$ and applying Lemma~\ref{lem:Stupin:caratheodoryCoeffDifference} gives
\begin{multline*}
    \operatorname{Re}\sum\limits_{j=1}^{p-1}\{h\}_j\{h\}_{p-j}
    =
    \frac12\sum\limits_{j=1}^{p-1}
    \bigl(|\{h\}_{p-j}|^2+|\{h\}_j|^2-|\{h\}_{p-j}-\overline{\{h\}_j}|^2\bigr)
    =\\=
    \sum\limits_{j=1}^{p-1}|\{h\}_j|^2
    -\frac12\sum\limits_{j=1}^{p-1}|\{h\}_{p-j}-\overline{\{h\}_j}|^2
    \geqslant
    \sum\limits_{j=1}^{p-1}|\{h\}_j|^2-2(p-1)(2-\operatorname{Re}\{h\}_p).
\end{multline*}
Substituting~\eqref{eq:Stupin:returnDeficit} here gives
\begin{equation}
    \label{eq:Stupin:returnQuadratic}
    \operatorname{Re}\sum\limits_{j=1}^{p-1}\{h\}_j\{h\}_{p-j}
    \geqslant\sum\limits_{j=1}^{p-1}|\{h\}_j|^2-O\left(\varepsilon\sum\limits_{j=1}^{p-1}|\{h\}_j|^2\right).
\end{equation}

Adding and subtracting $2-\varepsilon(2-\{h\}_p)$ inside the absolute value gives
\begin{multline}
    \label{eq:Stupin:returnModulusSquare}
    \left|-\frac{e^t}{t}\{f_{t,p,\varepsilon,h}\}_p\right|^2
    =\\=
    \left|2-\varepsilon(2-\{h\}_p)
    +\left(-\frac{e^t}{t}\{f_{t,p,\varepsilon,h}\}_p-\bigl(2-\varepsilon(2-\{h\}_p)\bigr)\right)\right|^2
    =\\=
    |2-\varepsilon(2-\{h\}_p)|^2
    +4\operatorname{Re}\left(-\frac{e^t}{t}\{f_{t,p,\varepsilon,h}\}_p-
    \bigl(2-\varepsilon(2-\{h\}_p)\bigr)\right)
    \\
    +2\operatorname{Re}\left(-\varepsilon(2-\overline{\{h\}_p})
    \left(-\frac{e^t}{t}\{f_{t,p,\varepsilon,h}\}_p-\bigl(2-\varepsilon(2-\{h\}_p)\bigr)\right)\right)
    \\
    +\left|-\frac{e^t}{t}\{f_{t,p,\varepsilon,h}\}_p-\bigl(2-\varepsilon(2-\{h\}_p)\bigr)\right|^2.
\end{multline}

Since $|2-\{h\}_p|\leqslant4$, it follows from~\eqref{eq:Stupin:returnDifference} that
\begin{multline}
    \label{eq:Stupin:returnProductBound}
    \left|2\operatorname{Re}\left(-\varepsilon(2-\overline{\{h\}_p})
    \left(-\frac{e^t}{t}\{f_{t,p,\varepsilon,h}\}_p-\bigl(2-\varepsilon(2-\{h\}_p)\bigr)\right)\right)\right|
    \leqslant\\\leqslant
    8\varepsilon\left|-\frac{e^t}{t}\{f_{t,p,\varepsilon,h}\}_p-\bigl(2-\varepsilon(2-\{h\}_p)\bigr)\right|
    =O\left(\varepsilon^3\sum\limits_{j=1}^{p-1}|\{h\}_j|^2\right).
\end{multline}
Moreover, since $\sum\limits_{j=1}^{p-1}|\{h\}_j|^2\leqslant4(p-1)$,
it follows from~\eqref{eq:Stupin:returnDifference} that
\begin{equation}
    \label{eq:Stupin:returnSquareBound}
    \left|-\frac{e^t}{t}\{f_{t,p,\varepsilon,h}\}_p-\bigl(2-\varepsilon(2-\{h\}_p)\bigr)\right|^2
    =
    O\left(\varepsilon^4\sum\limits_{j=1}^{p-1}|\{h\}_j|^2\right).
\end{equation}
Substituting~\eqref{eq:Stupin:returnProductBound} and~\eqref{eq:Stupin:returnSquareBound}
into~\eqref{eq:Stupin:returnModulusSquare} and using
\eqref{eq:Stupin:returnValue}, \eqref{eq:Stupin:returnExpansion}, and~\eqref{eq:Stupin:returnQuadratic}, we obtain
\begin{multline}
    \label{eq:Stupin:returnContradiction}
    4=\left|-\frac{e^t}{t}\{f_{t,p,\varepsilon,h}\}_p\right|^2
    =\\=
    |2-\varepsilon(2-\{h\}_p)|^2
    +4\operatorname{Re}\left(-\frac{e^t}{t}\{f_{t,p,\varepsilon,h}\}_p-
    \bigl(2-\varepsilon(2-\{h\}_p)\bigr)\right)
    +O\left(\varepsilon^3\sum\limits_{j=1}^{p-1}|\{h\}_j|^2\right)
    =\\=
    |2-\varepsilon(2-\{h\}_p)|^2
    -2t\varepsilon^2\operatorname{Re}\sum\limits_{j=1}^{p-1}\{h\}_j\{h\}_{p-j}
    +O\left(\varepsilon^3\sum\limits_{j=1}^{p-1}|\{h\}_j|^2\right)
    \leqslant\\\leqslant
    4-2t\varepsilon^2\sum\limits_{j=1}^{p-1}|\{h\}_j|^2
    +O\left(\varepsilon^3\sum\limits_{j=1}^{p-1}|\{h\}_j|^2\right).
\end{multline}

All constants in these estimates depend only on $t$ and $p$.
Thus there exists a number $\varepsilon_0\in(0,1/2)$, depending only on $t$ and $p$,
such that for $0<\varepsilon<\varepsilon_0$, \eqref{eq:Stupin:returnContradiction} implies
\[
    4\leqslant4-t\varepsilon^2\sum\limits_{j=1}^{p-1}|\{h\}_j|^2<4
    \quad\text{when}\quad
    \sum\limits_{j=1}^{p-1}|\{h\}_j|^2>0.
\]

If $\sum\limits_{j=1}^{p-1}|\{h\}_j|^2=0$, then~\eqref{eq:Stupin:returnDeficit}
gives $\operatorname{Re}\{h\}_p=2$, and $|\{h\}_p|\leqslant2$ implies $\{h\}_p=2$.
In this case, by Lemma~\ref{lem:Stupin:caratheodoryCoeffDifference}, we have
$\{h\}_{p-j}=\overline{\{h\}_j}$ for $1\leqslant j<p$.
Comparing the coefficients of $\varepsilon^2$ in~\eqref{eq:Stupin:returnExpansion} and~\eqref{eq:Stupin:secondVariation}
gives
\[
    \sum\limits_{j=1}^{p-1}|\{h\}_j|^2
    =
    \sum\limits_{j=1}^{p-1}\{h\}_j\{h\}_{p-j}
    =
    4\left(p\sum\limits_{k=1}^q\alpha_k^2-1\right)
    >0,
\]
where the strict inequality follows from~\eqref{eq:Stupin:secondVariationPositive} since $h\neq H_p$.
This contradicts $\sum\limits_{j=1}^{p-1}|\{h\}_j|^2=0$.

Thus~\eqref{eq:Stupin:returnValue} is impossible for $h\neq H_p$ and $0<\varepsilon<\varepsilon_0$.
Taking $n$ so that $1/n<\varepsilon_0$ contradicts the existence of
an $\varepsilon\in(0,1/n)$ satisfying~\eqref{eq:Stupin:returnValue}.
Consequently, $f_{t,p}$ gives a strict local maximum.

The rotation $z\mapsto e^{i\varphi/p}z$ preserves the class $C$ and maps basic neighborhoods to basic neighborhoods.
It commutes with $G_t$ and preserves the modulus of the $p$-th coefficient.
Consequently, the functions $F(e^{i\varphi}z^p,t)$, $\varphi\in[0,2\pi)$, give a strict local maximum
of the functional $I_{t,p}$ on $B_t$ in the topology $\tau_{B_t}$.

Since $(2te^{-t})'=2e^{-t}(1-t)$, we have $\max\limits_{t\geqslant0}{2te^{-t}}=2/e$,
with the maximum attained at $t=1$.
\end{proof}


\end{document}